\documentclass[11pt]{article}
\usepackage[letterpaper, margin=1in]{geometry}

\usepackage{amscd,amsfonts,amsmath,amssymb,amsthm,array,enumitem,fancyhdr,graphicx,multicol,scalerel,spectralsequences,subcaption,svg,tabularx,thmtools,thm-restate,tikz,tikz-cd}
\usepackage{xcolor}
\usepackage[colorlinks,
    linkcolor={violet},
    citecolor={olive!80!orange},
    urlcolor={blue!80!black}]{hyperref}
\usepackage[alphabetic]{amsrefs}

\theoremstyle{plain}
\newtheorem{theorem}{Theorem}[section]

\newtheorem{conjecture}[theorem]{Conjecture}
\newtheorem*{theorem*}{Theorem}

\theoremstyle{definition}
\newtheorem{definition}[theorem]{Definition}
\newtheorem*{definition*}{Definition}
\newtheorem{question}[theorem]{Question}

\theoremstyle{remark}

\newcommand{\Z}{\mathbb{Z}}
\newcommand{\Q}{\mathbb{Q}}

\newcommand{\Hom}{\operatorname{Hom}}

\newcommand{\ol}[1]{\overline{#1}}
\newcommand{\wh}[1]{\widehat{#1}}

\makeatletter
\newcommand{\colim@}[2]{%
  \vtop{\m@th\ialign{##\cr
    \hfil$#1\operator@font colim$\hfil\cr
    \noalign{\nointerlineskip\kern1.5\ex@}#2\cr
    \noalign{\nointerlineskip\kern-\ex@}\cr}}%
}
\newcommand{\colim}{%
  \mathop{\mathpalette\colim@{\rightarrowfill@\scriptscriptstyle}}\nmlimits@
}
\makeatother

\title{Triple linking and rational homology cobordism}
\author{Ryan Stees}
\date{}

\begin{document}

\maketitle


\begin{abstract}
If a rational homology 3-sphere $M$ bounds a rational homology 4-ball $W$, then the kernel of the inclusion-induced homomorphism $H_1(M;\Z)\to H_1(W;\Z)$ is a Lagrangian for the $\Q/\Z$-valued torsion linking form $\lambda_2$ on $H_1(M;\Z)$.
In this short paper, we prove that the Freedman-Krushkal triple torsion linking form $\lambda_3$ \cite{FreedmanKrushkal} vanishes on this Lagrangian under the assumption that $H_2(W;\Z)=0$.
We then pose several questions about topological rational homology cobordism.
\end{abstract}


\section{Introduction} \label{sec:intro}

A classical result states that if a rational homology 3-sphere $M$ bounds a topological rational homology 4-ball $W$, then the kernel $L$ of the inclusion-induced homomorphism $H_1(M;\Z)\to H_1(W;\Z)$ is a \emph{Lagrangian} for the torsion linking form \[\lambda_2:H_1(M;\Z)\times H_1(M;\Z)\to\Q/\Z.\]
In other words, $L\leq H_1(M;\Z)$ satisfies $L=L^\perp$ and $|L|=\sqrt{|H_1(M;\Z)|}$, where \[L^\perp:=\{[y]\in H_1(M;\Z)\,|\,\lambda_2([x],[y])=0 \text{ for all }[x]\in L\}.\]

Freedman-Krushkal recently developed a \emph{triple torsion linking form} $\lambda_3$ which is defined on triples $([x],[y],[z])$ of elements in $H_1(M;\Z)$, where $M$ is rational homology 3-sphere, such that $\lambda_2([x],[y])=\lambda_2([x],[z])=0\in\Q/\Z$. \cite{FreedmanKrushkal}.
In particular, for $M$ a rational homology sphere with $W$ and $L$ as above, $\lambda_3$ is defined on triples of elements in $L$.
In this paper, we consider the behavior of this form on $L$.
Our main result states that $\lambda_3$ vanishes on $L$ if $H_2(W;\Z)=0$.

\begin{restatable*}{theorem}{vanishing} \label{thm:vanishing}
Suppose the rational homology 3-sphere $M$ bounds a topological rational homology 4-ball $W$ with $H_2(W;\Z)=0$, and let $L\leq H_1(M;\Z)$ be the kernel of the inclusion-induced homomorphism $H_1(M;\Z)\to H_1(W;\Z)$. Then the triple torsion linking form $\lambda_3$ vanishes on $L$, that is, for any $[x],[y],[z]\in L$, we have $\lambda_3([x],[y],[z])=0\in\Q/\Z$.
\end{restatable*}

After recalling the Freedman-Krushkal triple torsion linking form in Section \ref{sec:form}, we prove Theorem \ref{thm:vanishing} and discuss an example in Section \ref{sec:vanishing}.
We conclude the paper in Section \ref{sec:open} by posing a number of questions which highlight aspects of topological rational homology cobordism that are not yet well-studied. \\


\noindent \textbf{Convention.} All homology groups will be assumed to have $\Z$ coefficients. \\


\noindent \emph{Acknowledgements.}
The author thanks Slava Krushkal for many helpful conversations.


\section{The triple torsion linking form} \label{sec:form}

In this section, we present the definition of the triple torsion linking form $\lambda_3$ from \cite{FreedmanKrushkal}.
Although this form is defined more generally, we will restrict to the case where our chosen triple of homology classes has pairwise vanishing classical linking in $\Q/\Z$.

Fix $t\in\Z$ which annihilates $H_1(M;\Z)$.
Let $[x],[y],[z]\in H_1(M;\Z)$ be three classes with pairwise vanishing linking, that is, \[\lambda_2([x],[y])=\lambda_2([x],[z])=\lambda_2([y],[z])=0\in\Q/\Z.\]
Represent $[x]$, $[y]$, and $[z]$ by oriented embedded curves $x$, $y$, and $z$ in $M$.
As $t\cdot[x]=0\in H_1(M;\Z)$, there is a map of a compact oriented surface $\Sigma$ into $M$ with $\partial\Sigma=tx$.
More precisely, $\partial\Sigma$ maps onto $x$ by a $t$-fold covering map and $\Sigma$ is otherwise embedded near $\partial\Sigma$.
As $\lambda_2([x],[y])=\lambda_2([x],[z])=0\in\Q/\Z$, we may assume $y$ and $z$ do not intersect $\Sigma$.

Define cohomology classes $\Phi,\Psi\in H^1(\Sigma)= H^1(\Sigma,\partial\Sigma)$ as follows: Define homomorphisms $\phi,\psi:H_1(\Sigma)\to\Z$ by $\phi([a])=A\cdot y$ and $\psi([a])=A\cdot z$, where $a\in C_1(M)$ is any representative of $[a]$ and $A\in C_2(M)$ is any 2-chain in $M$ with $\partial A=ta$.
These homomorphisms are independent of the choices of $a\in C_1(M)$ (hence defined on $H_1(M)$) and $A\in C_2(M)$ (hence well-defined).
They yield cohomology classes $\Phi$ and $\Psi$, respectively, in $H^1(\Sigma)=H^1(\Sigma,\partial\Sigma)$, with $\Phi\cup\Psi\in H^2(\Sigma,\partial\Sigma)=\Z$.

\begin{definition} \label{def:triple1}
The triple linking of the classes $[x]$, $[y]$, and $[z]$ is \[\lambda_3([x],[y],[z])=\frac{1}{t}(\Phi\cup\Psi)\in\Q/\Z.\]
\end{definition}

Freedman-Krushkal show that $\lambda_3$ is well-defined in $\Q/\Z$, independent of the many choices made in the definition.
Note that $\lambda_3$ does depend on the integer $t$ which we fixed earlier.
They also present a useful reformulation of Definition \ref{def:triple1}.
Suppose $\Sigma$ has genus $g$, and let $\{\gamma_i,\delta_i\}_{i=1}^g$ be a symplectic basis of oriented curves for $\Sigma$ such that each ordered pair $(\gamma_i,\delta_i)$ induces the given orientation of $\Sigma$.
As with the surface $\Sigma$ corresponding to the curve $x$, there are maps of compact oriented surfaces $C_i$ and $D_i$, $1\leq i\leq g$, into $M$ such that $\partial C_i=t\gamma_i$ and $\partial D_i=t\delta_i$ for all $i$.
As in \cite{FreedmanKrushkal}, we call $\Sigma$ the \emph{base surface} and the $C_i$ and $D_i$ \emph{second-stage surfaces}.

\begin{definition} \label{def:triple2}
The triple linking of the classes $[x]$, $[y]$, and $[z]$ is \[\lambda_3([x],[y],[z])=\frac{1}{t}\sum_{i=1}^g (C_i\cdot y)(D_i\cdot z)-(C_i\cdot z)(D_i\cdot y)\in\Q/\Z.\]
\end{definition}


\section{Vanishing theorem} \label{sec:vanishing}

A rational homology 3-sphere $M$ which embeds in $S^4$ decomposes the 4-sphere into two rational homology 4-balls $W_1$ and $W_2$.
A theorem of Hantzsche \cite{Hantzsche} states in part that the kernels of the two homomorphisms on $H_1$ induced by the inclusions $M\hookrightarrow W_1$ and $M\hookrightarrow W_2$ form \emph{dual Lagrangians} for the classical linking form $\lambda_2$ on $H_1(M)$.
Freedman-Krushkal show that the triple torsion linking form $\lambda_3$ also vanishes on each of these Lagrangians \cite{FreedmanKrushkal}.
In this section, we prove the main result of this paper, which extends the vanishing of $\lambda_3$ to the case where $M$ bounds any rational homology 4-ball $W$ with $H_2(W)=0$.

\vanishing

\begin{proof}
Our proof combines strategies seen in the proofs of the well-definedness of $\lambda_3$ and its vanishing on dual Lagrangians for $\lambda_2$ arising from embeddings in $S^4$ (Theorems 3.1 and 1.2, respectively, of \cite{FreedmanKrushkal}).

Suppose the rational homology 3-sphere $M$ bounds a topological rational homology 4-ball $W$ with $H_2(W)=0$, let $L=\ker(H_1(M)\to H_1(W))$, and let $t\in\Z$ be the integer annihilating $H_1(M)$ with respect to which we define $\lambda_3$.
As \[H_1(W,M)\cong H^2(W,M)\cong H_2(W)=0,\] where the first isomorphism follows from the Universal Coefficient Theorem and the second is Poincar\'{e} duality, the homology long exact sequence of the pair $(W,M)$ yields \[0\to H_2(W,M)\to H_1(M)\to H_1(W)\to 0.\]
In particular, $H_1(M)\to H_1(W)$ is surjective with kernel $L\cong H_2(W,M)$.
Thus, $t$ also annihilates $H_1(W)$.

Let $[x],[y],[z]\in L$, represented by oriented embedded curves $x$, $y$, and $z$.
By the classical result that $\lambda_2$ vanishes on $L$, the triple torsion linking $\lambda_3([x],[y],[z])$ is defined.
We will associate to $x$, $y$, and $z$ closed oriented surfaces $X$, $Y$, and $Z$ mapped into $W$ and define a quantity $\langle X,Y,Z\rangle_\Q\in\Z$.
We will then show that 
\begin{equation} \label{eq:lambda}
\lambda_3([x],[y],[z])=\frac{1}{t^3}\langle X,Y,Z\rangle_\Q\in\Q/\Z
\end{equation} and that 
\begin{equation} \label{eq:t3}
\langle X,Y,Z\rangle_\Q\in t^3\Z.
\end{equation}
Statements (\ref{eq:lambda}) and (\ref{eq:t3}) combine to show $\lambda_3([x],[y],[z])=0\in\Q/\Z$.

Denote the surface $\Sigma$ used in Definition~\ref{def:triple1} by $\Sigma_x$, and consider analogous compact oriented surfaces $\Sigma_y$ and $\Sigma_z$ mapped into $M$ with $\partial\Sigma_y=ty$ and $\partial\Sigma_z=tz$.
Denote the pushoff of $\Sigma_x$ into a collar neighborhood of $\partial W=M$ by $\Sigma_x^-\subset M\times\{-\varepsilon\}\subset M\times(-1,0]$, where $M\times\{0\}=\partial W$.
Because $[x],[y],[z]\in L$, the curves $x$, $y$, and $z$ bound compact oriented surfaces $S_x$, $S_y$, and $S_z$, respectively, mapped properly into $W$.
Let $S_x$, $S_y$, and $S_z$ have product structures on their intersections with $M\times(-1,0]$, and let $S_x^-$ be the surface $S_x\cap (W-(M\times(-\varepsilon,0])$, so that $\partial\Sigma_x^-=t\cdot (\partial S_x^-)$.
Consider the closed oriented surfaces $X=\Sigma_x^--tS_x^-$, $Y=\Sigma_y-tS_y$, and $Z=\Sigma_z-tS_z$ mapped into $W$.
As $H_2(W)=0$, the intersection number of $X$ with each of $Y$ and $Z$ vanishes.
We may therefore assume, after surgering the surfaces if necessary, that $X$ is disjoint from each of $Y$ and $Z$. 

As in Definition~\ref{def:triple2}, consider a symplectic basis $\{\gamma_i,\delta_i\}_{i=1}^g$ for $X$.
Because $t$ annihilates $H_1(W)$, there are maps of compact oriented (second-stage) surfaces $C_i$ and $D_i$ into $W$ such that $\partial C_i=t\gamma_i$ and $\partial D_i=t\delta_i$ for all $i$.
Define $\langle X,Y,Z\rangle_\Q$ by 
\begin{equation} \label{eq:rationalMatsumoto}
\langle X,Y,Z\rangle_\Q=\sum_{i=1}^g (C_i\cdot Y)(D_i\cdot Z)-(C_i\cdot Z)(D_i\cdot Y)\in\Z.
\end{equation}

This expression is denoted $\langle X,Y,Z\rangle_\Q^x$ in \cite{FreedmanKrushkal} and is one term of three in their ``rational Matsumoto pairing".
The expression $\langle X,Y,Z\rangle_\Q$ is defined on the \emph{chains} $X$, $Y$, and $Z$; we do not seek to define it on homology classes.
It is well-defined because there is a basis-independent formulation of $\langle X,Y,Z\rangle_\Q$ analogous to Definition \ref{def:triple1}.
We have \[\langle X,Y,Z\rangle_\Q=\ol{\Phi}\cup\ol{\Psi}\in \Z,\] where $\ol{\Phi}, \ol{\Psi}\in H^1(X)$ are analogous to $\Phi,\Psi\in H^1(\Sigma,\partial\Sigma)$ seen in Definition \ref{def:triple1}.

We now prove (\ref{eq:lambda}) holds.
As $\langle X,Y,Z\rangle_\Q$ is independent of the symplectic basis chosen for $X$, we may choose a basis $\{\gamma_i,\delta_i\}_{i=1}^g$ which splits into symplectic bases $\{\gamma_i,\delta_i\}_{i=1}^k$ and $\{\gamma_i,\delta_i\}_{i=k+1}^g$ for $\Sigma_x^-$ and the $t$ copies of $-S_x^-$, respectively.
Again by the well-definedness of $\langle X,Y,Z\rangle_\Q$, the second-stage surfaces corresponding to $\Sigma_x^-$ may be chosen to lie in $M\times\{-\varepsilon\}$.
A priori, intersections of the second-stage surfaces with $Y$ and $Z$ come in two types: those with $\Sigma_y$ and $\Sigma_z$, and those with $-S_y$ and $-S_z$ which come with multiplicity $t$.
By construction, the second-stage surfaces do not intersect $\Sigma_y$ and $\Sigma_z$.
We count intersections of the second-stage surfaces $C_i$ and $D_i$ with $-S_y$ and $-S_z$ separately for $1\leq i\leq k$ and $k+1\leq i\leq g$.

Consider the second-stage surfaces with $k+1\leq i\leq g$.
As $X$ consists of $t$ copies of the surface $-S_x^-$, each product of intersections of these second-stage surfaces with $-tS_y$ and $-tS_z$ appears as a summand of (\ref{eq:rationalMatsumoto}) $t$ times.
Each such term in (\ref{eq:rationalMatsumoto}) comes with an additional factor of $t^2$, one factor of $t$ coming from intersections with each of $-tS_y$ and $-tS_z$.
Thus, \[\sum_{i=k+1}^g(C_i\cdot Y)(D_i\cdot Z)-(C_i\cdot Z)(D_i\cdot Y)\in t^3\Z.\]

For $1\leq i\leq k$, the second-stage surfaces $C_i$ and $D_i$ lie in $M\times\{-\varepsilon\}$ and therefore intersect $-tS_y$ and $-tS_z$ in parallel copies of their boundaries.
More precisely, these intersections are counted as intersections with $-ty$ and $-tz$, respectively, in $M\times\{-\varepsilon\}$.
Thus, \[\sum_{i=1}^k(C_i\cdot Y)(D_i\cdot Z)-(C_i\cdot Z)(D_i\cdot Y)=t^2\sum_{i=1}^k(C_i\cdot y)(D_i\cdot z)-(C_i\cdot z)(D_i\cdot y),\] which implies \[\frac{1}{t^3}\sum_{i=1}^k(C_i\cdot Y)(D_i\cdot Z)-(C_i\cdot Z)(D_i\cdot Y)=\lambda_3([x],[y],[z])\in\Q/\Z\] by Definition \ref{def:triple2}.
Property (\ref{eq:lambda}) follows.

To prove (\ref{eq:t3}), we show first that we can modify the surface $X$ to make it nullhomologous in $E:=W-(\nu S_y\cup \nu S_z)$, where $\nu S_y$ and $\nu S_z$ are regular neighborhoods of the properly embedded surfaces $S_y$ and $S_z$.
Our modifications will change $\langle X,Y,Z\rangle_\Q$, which is additive in the first argument under connected sum, by an element of $t^3\Z$.
In particular, we will make $X$ nullhomologous in $E$ while preserving property (\ref{eq:lambda}).

We first study the Mayer-Vietoris sequence corresponding to the decomposition \[W=E\cup(\nu S_y\sqcup\nu S_z).\]
Because $H_3(W)\cong H^1(W,M)\cong\Hom(H_1(W,M),\Z)=0$ and $H_2(W)=0$ by assumption, this sequence yields the isomorphism \[ H_2(\partial(\nu S_y))\oplus H_2(\partial(\nu S_z))\xrightarrow{\cong} H_2(E).\]
Note that $H_2(\nu S_y)=H_2(\nu S_z)=0$ because $S_y$ and $S_z$ are surfaces with boundary.
Each of $H_2(\partial(\nu S_y))$ and $H_2(\partial(\nu S_z))$ is free abelian, generated by fundamental classes of tori $\alpha\times \mu$, where $\alpha$ represents an element of $H_1(S_y)$ (resp. $H_1(S_z)$) and $\mu$ is the meridian of $S_y$ (resp. $S_z$).
Thus, $H_2(E)$ is free abelian, generated by the images of these fundamental classes.

We arrange that $[X]=0\in H_2(E)$ by piping the surface $X$ into these tori $\alpha\times\mu$ with appropriate multiplicities.
We continue to call the modified surface $X$ and continue to denote its genus by $g$.
The effect of taking a connected sum with a single such torus $\alpha\times\mu$ increases the genus of $X$ by 1, and we expand the chosen symplectic basis for $X$ correspondingly using the curves $\alpha$ and $\mu$.
The resulting effect on $\langle X,Y,Z\rangle_\Q$ is seen in two new terms.

Suppose the torus $\alpha\times\mu$ represents a generator of $H_2(\partial(\nu S_y))$.
Let $A$ be a second-stage surface with $\partial A=t\alpha$, and take $t$ times the meridional disk $D^2$ as the second-stage surface for $\mu$.
As $D^2\cdot Z=0$, the first additional term $(A\cdot Y)(tD^2\cdot Z)$ in (\ref{eq:rationalMatsumoto}) vanishes.
For the second additional term $(A\cdot Z)(tD^2\cdot Y)$, recall that the second-stage surfaces do not intersect $\Sigma_y$ and $\Sigma_z$ by construction, so $A\cdot Z\in t\Z$.
The intersection $D^2\cdot S_y=1$ because $D^2$ is a meridional disk for $S_y$, so $(tD^2)\cdot Y=(tD^2)\cdot(-tS_y)=-t^2$.
Thus, $(A\cdot Z)(tD^2\cdot Y)\in t^3\Z$, and taking a connected sum with the torus $\alpha\times\mu$ changes $\langle X,Y,Z\rangle_\Q$ by an element of $t^3\Z$.
The argument for tori representing generators of $H_2(\partial(\nu S_z))$ is the same, so we have shown we can modify $X$ so that property (\ref{eq:lambda}) holds and $[X]=0\in H_2(E)$.

As $[X]=0\in H_2(E)$, $X$ bounds an oriented 3-manifold $N$ mapped into $E$.
We leverage this fact to prove (\ref{eq:t3}).
A standard application of Poincar\'{e} duality implies $K_\Q=\ker(H_1(X;\Q)\to H_1(N;\Q))$ is a (rational) Lagrangian for $H_1(X;\Q)$.
The intersection $K=K_\Q\cap H_1(X;\Z)$ is an integral Lagrangian, that is, a submodule $K\leq H_1(X;\Z)$ such that $H_1(X;\Z)/K\cong\Z^g$.
Represent the generators of $K$ by a \emph{geometric Lagrangian}, a collection of $g$ pairwise disjoint embedded curves $\alpha_i\subset X$ which form half of a symplectic basis for $X$.
As $[\alpha_i]\in K_\Q$ for $1\leq i\leq g$, for each $i$ there exists $n_i\neq 0\in\Z$ such that $n_i\alpha_i$ bounds a surface $A_i\subset N\subset E$.

Recall that $\langle X,Y,Z\rangle_\Q$ is independent of the symplectic basis chosen for $X$.
Choose such a symplectic basis $\{\gamma_i,\delta_i\}$ with $\gamma_i=\alpha_i$, and find $C_i$ with $\partial C_i=t\gamma_i$ as in the definition of $\langle X,Y,Z\rangle_\Q$.
As $H_2(W)=0$, $tA_i-n_iC_i$ is nullhomologous in $W$, hence \[n_i(C_i\cdot Y)=(n_iC_i)\cdot Y=(tA_i)\cdot Y=t(A_i\cdot Y)=0,\] where the final equality is a consequence of $A_i\subset E$.
Thus, $C_i\cdot Y=0$ for all $i$ and, by the same argument, $C_i\cdot Z=0$ for all $i$.
This proves (\ref{eq:t3}).
\end{proof}


\subsection{An example} \label{subsec:vanishing:ex}

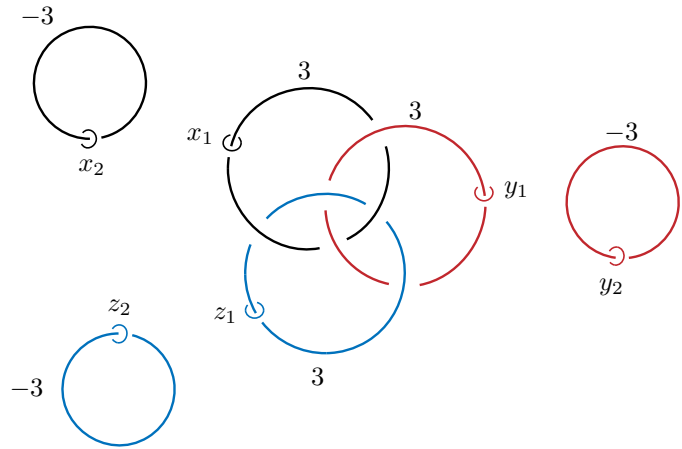
\begin{figure}[b]
\centering
\includegraphics[width=9cm]{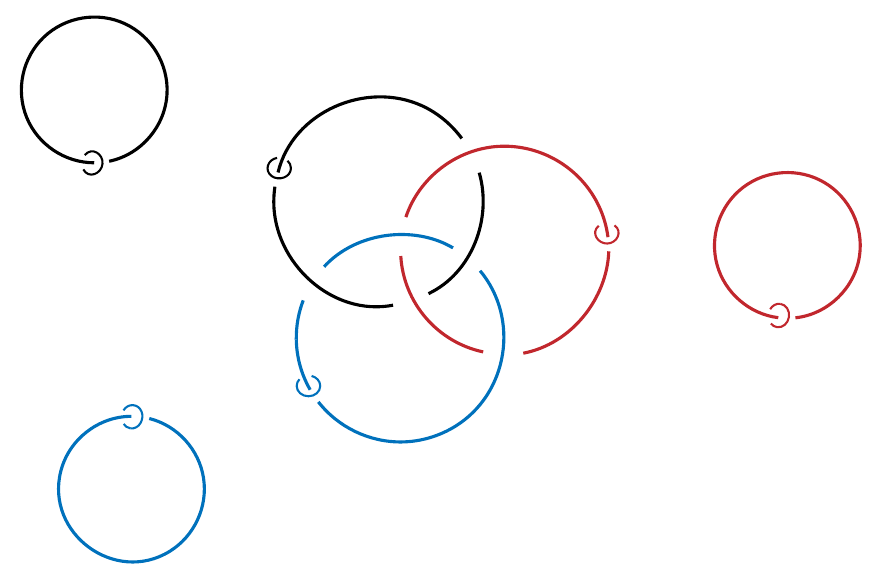}
\small
\put(-255,164){$-3$}
\put(-259,23){$-3$}
\put(-34,120){$-3$}
\put(-150,143){$3$}
\put(-108,128){$3$}
\put(-145,27){$3$}
\put(-233,108){$x_2$}
\put(-192,120){$x_1$}
\put(-36,63){$y_2$}
\put(-72,100){$y_1$}
\put(-222,55){$z_2$}
\put(-182,52){$z_1$}
\caption{A Kirby diagram for the manifold $M_{0}$, as seen in \cite{FreedmanKrushkal}, Figure 4.}
\label{fig:Example}
\end{figure}

Consider the 3-manifold $M_0$ from Section 7.1 of \cite{FreedmanKrushkal}, depicted in Figure \ref{fig:Example}.
The first homology of $M_0$ is isomorphic to $(\Z/3)^6$, generated by $\{[x_1],[y_1],[z_1],[x_2],[y_2],[z_2]\}$, the classes of the meridians to the surgery curves .
The torsion linking form of $M_0$ is isomorphic to that of $\#^3(L(3,1)\#-L(3,1))$.
Let $x=x_1+x_2$, $y=y_1+y_2$, and $z=z_1+z_2$.
The classes $[x]$, $[y]$, and $[z]$ generate a Lagrangian $L$ for $\lambda_2$.
Freedman-Krushkal compute that $\lambda_3([x],[y],[z])=\frac{1}{3}\in\Q/\Z$.
Thus, by Theorem \ref{thm:vanishing}, $M_0$ cannot bound a rational homology ball $W$ with $H_2(W)=0$ and $L=\ker(H_1(M_0)\to H_1(W))$.

On the other hand, Freedman-Krushkal exhibit Lagrangians for $\lambda_2$ on which $\lambda_3$ \emph{does} vanish.
For example, $\lambda_3$ vanishes on the Lagrangian generated by $\{[l_1],[l_2],[l_3]\}$, where
\begin{align*}
l_1&=x_2+y_2+z_2 \\
l_2&=y_1-z_1+y_2-z_2 \\
l_3&=x_1+y_1+z_1.
\end{align*}
We are therefore unable to use Theorem \ref{thm:vanishing} to obstruct $M_0$ from bounding any rational homology ball $W$ with $H_2(W)=0$.
In particular, we are unable to obstruct $M_0$ from bounding a \emph{ribbon} $\Q$-homology ball \cite{DLVSW}.

One might hope to construct a 3-manifold related to $M_0$ for which $\lambda_3$ does not vanish on \emph{any} Lagrangian for $\lambda_2$.
More specifically, consider all framed links which are the result of combinations of clasper surgeries on triples of surgery curves for $\#^3(L(3,1)\#-L(3,1))$.
Performing surgery on each of these framed links yields a family of 3-manifolds $\{M_v\}$, $v\in(\Z/3)^{20}$.
Freedman-Krushkal compute, however, that each of the $3^{20}\approx 3.5\text{ billion}$ resulting 3-manifolds $M_v$ has a Lagrangian on which $\lambda_3$ vanishes \cite{FreedmanKrushkal}.
It turns out that the same can be said for the analogous 3-manifolds constructed from $\#^3(L(2,1)\#-L(2,1))$.


\section{Open questions} \label{sec:open}

We conclude this paper with several open questions arising from the mysteriously subtle behavior of $\lambda_3$.
To the author's knowledge, the following have not been considered in the literature.


\subsubsection*{Surfaces in rational homology balls}

The condition $H_2(W)=0$ in Theorem \ref{thm:vanishing} is used to conclude that the integer $t\in\Z$ which annihilates $H_1(M)$ also annihilates $H_1(W)$.
This assumption is stronger than necessary; one can adapt the proof of Theorem \ref{thm:vanishing} in the case where this conclusion holds even if $H_2(W)\neq 0$.
In fact, one does not need this conclusion globally; it is enough to know that curves representing classes in $L\leq H_1(M)$ bound surfaces in $W$ whose homology is in the image of $H_1(M)\to H_1(W)$.
\begin{question} \label{question:deep}
Given a rational homology 3-sphere $M$ and a bounding rational homology 4-ball $W$, under what circumstances does a representative of a given element of $\ker(H_1(M)\to H_1(W))$ bound a surface in $W$ whose first homology is in the image of $H_1(M)\to H_1(W)$?
Classes in this kernel which do not satisfy this condition might be called \emph{deeply nullhomologous} in $W$.
\end{question}


\subsubsection*{Bounding Lagrangians}

The example from Section \ref{subsec:vanishing:ex} exhibits a rational homology sphere $M_0$ and a Lagrangian $L$ for $\lambda_2^{M_0}$ such that $M_0$ cannot bound a rational homology ball $W$ with $H_2(W)=0$ and $L=\ker(H_1(M_0)\to H_1(W))$.
This example shows that $\lambda_3$ places algebraic restrictions on the kinds of rational homology balls which rational homology spheres bound and the ways in which they can bound them.

\begin{question} \label{question:bounding}
Suppose $M$ is a rational homology 3-sphere which bounds \emph{some} (topological) rational homology 4-ball.
For each Lagrangian $L\leq H_1(M;\Z)$ for $\lambda_2^M$, does $M$ necessarily bound a rational homology 4-ball $W_L$ with $L=\ker(H_1(M)\to H_1(W_L))$?
If not, for which Lagrangians is this the case?
\end{question}


\subsubsection*{The topological rational homology cobordism group}

The behavior of $\lambda_2^M$ when $M$ bounds a rational homology 4-ball enables the definition of a homomorphism\[\Lambda:\wh{\Theta}_3^\Q\to W(\Q/\Z)\] from the 3-dimensional topological rational homology cobordism group $\wh{\Theta}_3^\Q$ of rational homology 3-spheres to the Witt group $W(\Q/\Z)$ which sends the class $[M]$ of a 3-manifold $M$ to the class $[\lambda_2^M]$ of its torsion linking form.

We will call the group $W(\Q/\Z)$ the \emph{algebraic rational homology cobordism group}, and a rational homology sphere $M$ with $[M]\in\ker\Lambda$ \emph{algebraically nullbordant}.
This builds an analogy to the topological concordance group of knots in $S^3$ and its corresponding \emph{algebraic concordance group}, the group of Witt classes of Seifert forms.
Indeed, the branched double cover of $S^3$ branched along an algebraically slice knot is an algebraically nullbordant homology 3-sphere. 

Work of Kawauchi-Kojima \cite{KawauchiKojima} implies the homomorphism $\Lambda$ is surjective. 
Kim-Livingston comment that $\Lambda$ is conjecturally an isomorphism \cite{KimLivingston}.

\begin{conjecture}[\cite{KimLivingston}] \label{conj:isomorphism}
If $M$ is algebraically nullbordant, then $M$ is nullbordant in $\wh{\Theta}^\Q_3$.
\end{conjecture}

A positive answer to the following question, along with an understanding of Question \ref{question:deep}, could lead to a counterexample.
\begin{question} \label{question:counterex}
Does there exist a rational homology sphere $M$ such that $\lambda_3$ is nonvanishing on \emph{each} Lagrangian for $\lambda_2$?
\end{question}


\bibliography{references.bib}


\vspace{1em}

\noindent
Ryan Stees, {University of Virginia}\\
Email: \texttt{\href{mailto:rs2sf@virginia.edu}{rs2sf@virginia.edu}}\\
URL: \url{https://www.ryanstees.com}

\end{document}